\documentclass{amsart}
\usepackage{amssymb}

\newtheoremstyle{notauto}{}{}{\itshape}{}{\bfseries}{.}{0.5em}{\thmnote{#3}}
\theoremstyle{notauto}

\theoremstyle{definition}

\theoremstyle{remark}

\begin{document}
\title{On abelian group actions with TNI-centralizers}
\author{G\"{u}l\.{I}n Ercan$^{*}$}
\address{G\"{u}l\.{I}n Ercan, Department of Mathematics, Middle East
echnical University, Ankara, Turkey}
\email{ercan@metu.edu.tr}
\thanks{$^{*}$Corresponding author}
\thanks{This work has been supported by the Research Project T\"{U}B\.{I}TAK
114F223.}
\author{\.{I}sma\.{I}l \c{S}. G\"{u}lo\u{g}lu}
\address{\.{I}sma\.{I}l \c{S}. G\"{u}lo\u{g}lu, Department of Mathematics, Do%
\u{g}u\c{s} University, Istanbul, Turkey}
\email{iguloglu@dogus.edu.tr}
\maketitle

\begin{abstract}  A subgroup $H$ of a group $G$ is said to be a TNI-subgroup if $N_{G}(H)\cap
H^g=1$ for any $g\in G\,\backslash \,N_{G}(H).$ Let $A$ be an abelian group acting coprimely on the finite group $G$ by automorphisms in such a way that $C_G(A)=\{g\in G : g^a=g $\, for all $a\in A\}$ is a solvable TNI-subgroup of $G$. We prove that $G$ is a solvable group with Fitting length $h(G)$ is at most $h(C_G(A))+\ell(A)$. In particular $h(G)\leq \ell(A)+3$ whenever $C_G(A)$ is nonnormal. Here, $h(G)$ is the Fitting length of $G$ and $\ell(A)$ is the number of primes dividing $A$ counted with multiplicities. 
\end{abstract}

\section{introduction}  Throughout the paper all groups are finite and $h(G)$ denotes the Fitting length of the group $G$. A subgroup $H$ of a group $G$ is said to be a TNI-subgroup if $N_{G}(H)\cap H^g=1$ for any $g\in G\,\backslash \,N_{G}(H).$ In particular, every normal subgroup is a $TNI$-subgroup. In \cite{EG}, we studied the consequences of the action of a group $A$ on the group $G$ in case where $C_G(A)$ is a $TNI$-subgroup, and obtained the following two results:\\

\textbf{Theorem A.} \textit{ Let $A$ be a group that acts coprimely on the group $G$ by automorphisms. If $C_G(A)$ is a solvable $TNI$-subgroup of $G$, then $G$ is solvable.}\\

\textbf{Theorem B.} \textit{Let $A$ be a coprime automorphism of prime order of a finite solvable group $G$ such that $C_G(A)$ is a $TNI$-subgroup. Then $h(G)\leq h(C_G(A))+1.$ In particular, $h(G)\leq 4$ when $C_G(A)$ is nonnormal.}\\

In the present paper we extend Theorem B to the case where $A$ is abelian, namely we prove\\

\textbf{Theorem.} \,\,\textit{ Let $A$ be an abelian group acting coprimely on the finite group $G$ in such away that $C_G(A)$ is a solvable TNI-subgroup of $G$. Then $G$ is a solvable group with $h(G)\leq h(C_G(A))+\ell(A)$ where $\ell(A)$ is the number of primes dividing $A$ counted with multiplicities. In particular $h(G)\leq \ell(A)+3$ whenever $C_G(A)$ is nonnormal.}\\

This is achieved by applying the same techniques used in \cite{Tur1} in order to prove that $h(G)\leq \ell(A)$ if $A$ acts coprimely and fixed point freely on the group $G$ and for every proper subgroup $D$ and every $D$-invariant section $S$ of $G$ such that $D$ acts irreducibly on $S$, there is $v\in S$ with $C_D(S)=C_D(v)$, that is, $A$ acts with regular orbitd on $G$. In the light of this result we asked whether the conclusion of Theorem B is true if $A$ is not necessarily abelian, but acts with regular orbits on $G$. The main difficulty forcing us to study under the assumption that $A$ is abelian arises from the fact that a homomorphic image of a $TNI$-subgroup need not be a $TNI$-subgroup. 

\section{proof of the theorem} 

The group $G$ is solvable by Theorem A in \cite{EG}. It remains to show that $h(G)\leq h(C_G(A))+\ell(A).$  Suppose false and let\, $G,A$  be a counterexample with $|G|$ minimum.\\

Suppose that $C_G(A)$ is normal in $G$. Then the fixed point free action of $A$ on the group $G/C_G(A)$ yields that $h(G/C_G(A))\leq \ell(A)$ by the main theorem of \cite{Tur1}. So $h(G)\leq h(C_G(A))+\ell(A)$, a contradiction. Therefore by Theorem 2.2 in \cite{EG} we may assume that\\

\textit{$(1)$ $C_G(A)$ is a nonnormal subgroup of $G$ acting Frobeniusly on a section $M/N$ of $G$.}\\

We also have\\

\textit{$(2)$ There is an $A$-tower $\hat{P}_i, i=1,\ldots ,t$ where $t=h(G)$ satisfying the following conditions (see \cite{Tur2})}:

\vspace{1mm} \textit{(a)} $\hat{P}_{i}$\textit{\ is an $A$-invariant }$p_{i}
$\textit{-subgroup, }$p_{i}$\textit{\ is a prime, }$p_{i}\neq p_{i+1},$%
\newline
\textit{\ for }$i=1,\ldots ,t-1$\textit{;}

\textit{(b)} $\hat{P}_{i}\leq N_{G}(\hat{P}_{j})$\textit{\ whenever }$i\leq j
$\textit{;}

\textit{(c)} $P_{t}=\hat{P}_{t}$\textit{\ and }$P_{i}=\hat{P}_{i}/C_{\hat{P}%
_{i}}(P_{i+1})$ \textit{\ for }$i=1,\ldots ,t-1$ \newline
\textit{\ and }$P_{i}\neq1$ \textit{\ for }$i=1,\ldots ,t$\textit{;}

\textit{(d)} $\Phi(\Phi(P_{i}))=1$\textit{, }$\Phi(P_{i})\leq Z(P_{i})$%
\textit{, and exp}$(P_{i})=p_{i}$\textit{\ when }$p_{i}$\textit{\ is odd } 
\newline
\textit{\ for} $i=1,\ldots ,t$ \textit{;}

\textit{(e)} $[\Phi(P_{i}),\hat{P}_{i-1}]=1$\textit{\ and }$[P_{i},P_{i-1}]=P_{i}
$\textit{\ for }$i=1,\ldots ,t$\textit{;}

\textit{{(f)} If $S\leq \hat{P}_{i}$ for some $i$, $S$ is normalized by $\hat{P}_{i-1}\ldots \hat{P}_{1}A$ and its image in ${P}_{i}$ is not contained in $\Phi(P_{i})$, then $S=\hat{P}_{i}$.}\\

By Lemma 2.1 in \cite{EG} the group $\prod_{i=1}^{t}{\hat{P}_{i}}$ is of Fitting length $t$ and it satisfies the hypothesis of the theorem. It follows now by induction that 
\\

\textit{$(3)$ $G=\prod_{i=1}^{t}{\hat{P}_{i}}$.}\\

Suppose that $C_{\hat{P}_{t}}(A)\ne 1$. Then we have $M/N=[M,C_{\hat{P}_{t}}(A)]N/N$ due to the Frobenius action of $C_{\hat{P}_{t}}(A)$ on $M/N.$ It follows that $M/N\leq \hat{P}_{t}N/N\cap M/N=1$ as $\hat{P}_{t}\lhd G$ and $p_{t}$ is coprime to $|M/N|$. This contradiction shows that\\

\textit{$(4)$  $C_{\hat{P}_{t}}(A)=1$.}\\

Set $H=\prod_{i=1}^{t-1}{\hat{P}_{i}}$.  Pick a nontrivial subgroup $C$ of $A$. Set $S=[\hat{P}_{t-1},C]^{H}$. Clearly $S\leq \hat{P}_{t-1}$ is normalized by $\hat{P}_{t-1}\ldots \hat{P}_{1}A$. Now the image of $S$ in $P_{t-1}$ is $[{P}_{t-1},C]^{H}$. Suppose that $[{P}_{t-1},C]^{H}$ is contained in $\Phi({P}_{t-1})$. It follows that $[{P}_{t-1},C]\leq \Phi({P}_{t-1})$ and so $[{P}_{t-1},C ]=1$ due to coprimeness. By the three subgroup lemma  $[{P}_{t-2},C,{P}_{t-1}]=1$ whence $[{P}_{t-2},C]=1$. Repeating the same argument one gets $[{P}_{i},C]=1$ for each $i< t$. Now the group $X=\prod_{i=1}^{t-1}{C_{\hat{P}_{i}}}(C)$ is of Fitting length $t-1$ on which $A/C$ acts in such a way that $C_X(A/C)$ is a TNI-subgroup. By induction we get $t-1\leq f(C_X(A))+\ell(A/C)$. It then follows that $t\leq f(C_G(A))+\ell(A)$. This contradiction shows that $[{P}_{t-1},C]^{H}$ is not contained in $\Phi({P}_{t-1})$. By $(2)$ part $f$ we have $S=\hat{P}_{t-1}$. This shows that $[\hat{P}_{t-1},C]^{\hat{P}_{t-1}\ldots \hat{P}_{1}}={\hat{P}_{t-1}}$ for every subgroup $C$ of $A$ with $\ell(C)\geq 1$.\\

Next let $D\leq A$ with $\ell(D)\geq 2$ and $Y=\prod_{i=1}^{t-2}{\hat{P}_{i}}$. Set $T=[\hat{P}_{t-2},D]^{Y}$. Clearly $T$ is $YA$-invariant. If the image $[{P}_{t-2},C]^{Y}$ of $T$ in $P_{t-2}$ is contained in $\Phi({P}_{t-2})$, then we can show by an argument similar as in the above paragraph that $[{P}_{i},D]=1$ for each $i< t-1$. Then $Z=\prod_{i=1}^{t-2}{C_{\hat{P}_{i}}}(D)$ is a group of Fitting length $t-2$ on which $A/D$ acts in such a way that $C_Z(A/D)$ is a TNI-subgroup. It follows by induction that  $t\leq f(C_G(A)+\ell(A)$, which is not the case. Therefore $T= \hat{P}_{t-2}$ by $(2)$ part $(f)$. Thus we have \\

\textit{$(5)$ $[\hat{P}_{t-1},C]^{\hat{P}_{t-1}\ldots \hat{P}_{1}}={\hat{P}_{t-1}}$ for every subgroup $C$ of $A$ with $\ell(C)\geq 1$ and $[\hat{P}_{t-2},D]^{\hat{P}_{t-2}\ldots \hat{P}_{1}}= \hat{P}_{t-2}$ for every subgroup $D$ of $A$ with $\ell(D)\geq 2$. }\\

Let now $S$ be an $H$-homogeneous component of the irreducible $HA$- module $V={P}_{t}/\Phi({P}_{t})$. Notice that $\hat{P}_{t-1}$ acts nontrivially on each $H$-homogeneous component of $V.$ Set $B=N_A(S)$. Then $S$ is an irreducible $HB$-module such that $C_S(B)=0$. By the Fong-Swan theorem\\

\textit{$(6)$ We may take an irreducible $\mathbb{C}HB$-module $M$ such that $M|_H$ is homogeneous, $Ker_H(M)=Ker_H(S)$ and $C_{M}(B)=0.$ Among all pairs $(M_{\alpha}, C_{\alpha})$ such that $1\ne C_{\alpha}\leq B$, $M_{\alpha}$ is an irreducible  $HC_{\alpha}$ submodule of $M|_{HC_{\alpha}}$ and $C_ {M_{\alpha}}(C_{\alpha})=0$, choose $(M_1,C)$ with $|C|$ minimum. Then $C_{M_1}(C_0)\ne 0$ for $1\ne C_0 < C$ and 
$Ker_H(M)= Ker_H(M_1).$ }\\

Set $\bar{H}=H/Ker_H(M).$ Suppose that $\bar{U}\ne 1$ is an abelian subgroup of $\bar{H}$ and that $\bar{U}\lhd \bar{H}C.$ Let $U$ be the preimage in $H$ of $\bar{U}.$ Since $M_1|_H$ is homogeneous, by Glauberman's lemma there is a homogeneous component $N_1$ of $M_1|_U$ such that $C\leq N_{HC}(N_1)$. Set $H_1=N_H(N_1).$ Then we have $H_1C=N_{HC}(N_1)$. Now $[U,C]\leq Ker_H(N_1)$.  
By Proposition 4.1 in \cite{Tur4}, $H=({\bigcap_{x\in HC}{H_1}^x})C_H(C)$. It follows that $M_1=\Sigma (N_1)^{x}$ for $x\in C_H(C).$ Notice that $[U,C]\leq Ker_H(N_1).$ Thus $[U,C]=[U,C]^{x}\leq Ker_H(N_1)^{x}$ and so $[U,C]\leq Ker_H(M_1)$. Then we have\\

\textit{$(7)$ If $\bar{U}$ is an abelian subgroup of $\bar{H}$ such that $\bar{U}\lhd \bar{H}C$ where $\bar{H}=H/Ker_H(M)$ then $[\bar{U},C]=1.$}\\

By $(5)$ we have $[\hat{P}_{t-1},C]^{H}C_{\hat{P}_{t-1}}({P}_{t})=\hat{P}_{t-1}$ and hence $[\hat{P}_{t-1},C]\not \leq Ker(M)$. Therefore by $(7)$, $\hat{P}_{t-1}Ker_H(M)/Ker_H(M)$ is nonabelian. As $P_1$ is elementary abelian we have $t>2.$ Set $\Phi/Ker_H(M)=\Phi(\hat{P}_{t-1}Ker_H(M)/Ker_H(M))$ and let $H_1=C_{\hat{P}_{t-2}\ldots \hat{P}_{1}}(\bar{\Phi})$. By $(7)$, $[\bar{\Phi},C]=1.$  
Now $H_1C=C_{\hat{P}_{t-2}\ldots \hat{P}_{1}C}(\bar{\Phi})\lhd \hat{P}_{t-2}\ldots \hat{P}_{1}C.$ Hence $[\hat{P}_{t-2}\ldots \hat{P}_{1},C]\leq H_1$. By the coprimeness we have \\

\textit{$(8)$ $t>2$ and $H_1C_{\hat{P}_{t-2}\ldots \hat{P}_{1}}(C)=\hat{P}_{t-2}\ldots \hat{P}_{1}$ where $\bar{\Phi}=\Phi(\hat{P}_{t-1}Ker(M)/Ker(M))$ and $H_1=C_{\hat{P}_{t-2}\ldots \hat{P}_{1}}(\bar{\Phi})$.}\\

We have $\hat{P}_{t-2}\leq H_1$ by $(2)$ part $(e)$. Set $Q=\hat{P}_{t-2}$. Let $D\leq C$ be such that $\ell(D)\geq 2.$ Then by $(8)$, $[Q,D]^{H_1}=[Q,D]^{\hat{P}_{t-2}\ldots \hat{P}_{1}}.$ By $(5)$ we have\\

\textit{$(9)$ If $D\leq C$ with $\ell(D)\geq 2$ then $Q=[Q,D]^{H_1}$ where $Q=\hat{P}_{t-2}$.}\\

Let $N$ be a homogeneous component of $M_1|_S$. Then $N$ is normalized by $\hat{P}_{t-1}H_1C$ since $H_1C=C_{\hat{P}_{t-2}\ldots \hat{P}_{1}C}(\bar{\Phi}).$ Set $P_0=\hat{P}_{t-1}/\hat{P}_{t-1}\cap Ker_H(N)$. Then $N$ is a $P_0H_1C$-module. Notice that $H_1C$ centralizes $\Phi(P_0)$ and hence $N|_{\Phi(P_0)}$ is homogeneous. Then $\Phi(P_0)$
is cyclic. We also have $\Phi(P_0)$ is elementary abelian by $(2)$ part $(d)$ and hence $|\Phi(P_0)|\leq p.$ Recall that ${(\hat{P}_{t-1})'}\not \leq Ker_{H}(M_1)$ and hence $P_0$ is nonabelian. Thus we have ${P_0}'=\Phi(P_0)$. Note also that $P_0/Z(P_0)$ is elementary abelian. As $P_{t-1}/\Phi(P_{t-1})$ is irreducible $\hat{P}_{t-2}\ldots \hat{P}_{1}A$-module, it is completely reducible as $H_1$-module because $H_1$ is subnormal in $\hat{P}_{t-2}\ldots \hat{P}_{1}A$. It follows that $P_0/\Phi(P_0)$ is $H_1$-completely reducible. It follows by Maschke's theorem that it is $H_1C$-completely reducible. 

Suppose that $1\ne D\leq C$. Then $[P_0,D]\leq[P_0,D]^{P_0H_1}\lhd P_0H_1C_{\hat{P}_{t-2}\ldots \hat{P}_{1}}(C).$ By $(8)$, $P_0H_1C_{\hat{P}_{t-2}\ldots \hat{P}_{1}}(C)=P_0\hat{P}_{t-2}\ldots \hat{P}_{1}$ and hence by $(5)$ we get $P_0=[P_0,D]^{H_1}$. Now we apply Theorem 1.1 in \cite{Tur3} by letting $G=P_0H_1$, $A=C$, $P=P_0$, $Q=\hat{P}_{t-2}$ and $\chi$ as the character afforded by $N$. This leads to  $C_N(C)\ne 0,$ which is a contradiction completing the proof.

\end{document}